\documentclass[oupthm]{CUP-JNL-BCM}%

\usepackage{graphicx}
\usepackage{multicol}
\usepackage{amsmath,amssymb,amsfonts, mathtools}
\usepackage{mathrsfs}
\usepackage{rotating}
\usepackage{appendix}
\usepackage[numbers]{natbib}
\usepackage{changepage,afterpage}
\usepackage{etruscan}
\usepackage{tikz, siunitx, newunicodechar}
\usepackage{booktabs,multirow}
\usepackage{bm,bbm}%
\usepackage{braket}

\usepackage[utf8]{inputenc}

\usepackage{array}
\usepackage{xcolor}
\usepackage{colortbl}

\allowdisplaybreaks[4]

\theoremstyle{oupplain}
\newtheorem{theorem}{Theorem}[section]

\newtheorem{corollary}[theorem]{Corollary}
\theoremstyle{oupdefinition}

\theoremstyle{oupremark}
\newtheorem{remark}[theorem]{Remark}

\theoremstyle{oupproof}

\numberwithin{equation}{section}

\articletype{RESEARCH ARTICLE}

\begin{document}
\setlength{\abovedisplayskip}{3pt}
\setlength{\belowdisplayskip}{3pt}

\begin{Frontmatter}
\title{Majorization in some symplectic weak supermajorizations}
\author{Shaowu Huang}
\address{
\orgname{School of Mathematics and Finance, Putian University},
 \orgaddress{\street{Putian Fujian}, 
 \state{People's Republic of China},
  \postcode{351100}}
  \email{shaowu2050@126.com}}
\author{Hemant K. Mishra}
\address{
\orgname{School of Electrical and Computer Engineering, Cornell University}, 
\orgaddress{\street{Ithaca}, \state{New York}, \postcode{14850}}
\email{hemant.mishra@cornell.edu}}
\keywords[2020 Mathematics Subject Classification]{15B48, 15A18}
\keywords{Symplectic matrix, symplectic eigenvalue, Williamson's theorem, majorization, weak supermajorization}
\abstract{Symplectic eigenvalues are known to satisfy analogs of several classic eigenvalue inequalities.
    Of these is a set of  weak supermajorization relations concerning symplectic eigenvalues that are weaker analogs of some majorization relations corresponding to eigenvalues.
    The aim of this letter is to establish necessary and sufficient conditions for the saturation of the symplectic weak supermajorization relations by  majorization.}
\end{Frontmatter}
\section{Introduction}
    A $2n \times 2n$ real matrix $M$ is said to be symplectic if it satisfies $M^TJ_{2n}M = J_{2n}$, where 
        $J_{2n} \coloneqq
            \begin{psmallmatrix}
             0 & I_n \\
             -I_n & 0
            \end{psmallmatrix}
        $, 
    $I_n$ being the identity matrix of size $n$.
    Williamson's theorem \cite{williamson1936algebraic} states that for every $2n \times 2n$ real symmetric positive definite matrix $A$, there exists a symplectic matrix $M$ such that
    \begin{align}\label{eq:Williamson_decomposition}
        M^T A M 
            &=
            \begin{pmatrix}
             D & 0 \\
             0 & D
            \end{pmatrix},
    \end{align}
    where $D$ is an $n \times n$ diagonal matrix with positive diagonal entries, called the symplectic eigenvalues of $A$.
    See \cite{folland1989harmonic, simon1999congruences, ikramov2018symplectic} for elementary proofs of Williamson's theorem. 
    
    Symplectic eigenvalues are known to satisfy analogs of several classic eigenvalue inequalities \cite{bhatia2015symplectic, HIAI2018129, mishra2020first, bhatia2020schur, bhatia_jain_2021, jain2021sums, jm, paradan2022, mishra2023, sags_2021, huang2023}.
    Of these is a set of weak supermajorization relations concerning symplectic eigenvalues that are weaker analogs of some majorization relations corresponding to eigenvalues.
    The aim of this letter is to establish necessary and sufficient conditions for the saturation of the symplectic  weak supermajorization relations by majorization.

    Let $\operatorname{Pd}(2n)$ denote the set of $2n \times 2n$ real symmetric positive definite matrices.
    We denote by $\operatorname{Sp}(2n)$ the set of $2n \times 2n$ real symplectic matrices.
    The set $\operatorname{Sp}(2n)$ is a group under matrix multiplication, known as the  symplectic group.
    The symplectic group is a non-compact subset of the special linear group and is closed under transpose \cite{dms}.
    We call $M \in \operatorname{Sp}(2n)$ an {\it orthosymplectic} matrix if it is an orthogonal matrix, i.e., $M^TM=I_{2n}$.
    For $A \in \operatorname{Pd}(2n)$, let $\operatorname{Sp}(2n; A) \subset \operatorname{Sp}(2n)$ denote the set of symplectic matrices diagonalizing $A$ in the sense of Williamson's theorem \eqref{eq:Williamson_decomposition}. 
    We denote by $d_1(A), \ldots, d_n(A)$ the symplectic eigenvalues of $A$ arranged in non-decreasing order, which are uniquely determined.
    We say that $A$ is orthosymplectically diagonalizable in the sense of Williamson's theorem if there exists an orthosymplectic matrix $M$ diagonalizing $A$ as in \eqref{eq:Williamson_decomposition}.

    Given any vector $x$ in $\mathbb{R}^n$,  we shall denote by $\operatorname{diag}(x)$ the $n \times n$ diagonal matrix whose diagonal entries are given by the entries of $x$.
    Denote the entries of $x$ in the ascending order by $x^{\uparrow}_1 \leq  \cdots \leq x^{\uparrow}_n$. 
    Let $x,y$ be two vectors in $\mathbb{R}^n$. We say that $x$ is weakly supermajorized by $y$, written as $x \prec^w y$, if 
    \begin{align}\label{eq:supmajorization}
    \sum_{i=1}^k x^{\uparrow}_i \geq \sum_{i=1}^k y^{\uparrow}_i, \quad \text{for } k=1,\ldots, n.
    \end{align}
    In addition, if \eqref{eq:supmajorization} is an equality for $k=n$,  then $x$ is said to be majorized by $y$ and is written as $x \prec y$. 
    An $n \times n$ real matrix $E$ is said to be a doubly stochastic matrix if its $(i,j)$th entries $E_{ij}$ are non-negative for $i,j=1,\ldots, n$ such that
    \begin{align}
        \sum_{j=1}^n E_{ij}=1, \quad \text{for } i=1, \ldots, n, \\
        \sum_{i=1}^n E_{ij}=1, \quad \text{for } j=1,\ldots, n.
    \end{align}
    It is well-known from the theory of majorization that $x \prec y$ if and only if there exists an $n \times n$ doubly stochastic matrix $E$ such that $x=Ey$.
    See \cite[Theorem~1.3]{ando1989}.

\subsection{Symplectic Schur--Horn weak supermajorizations}
    Several weak supermajorization relations between the entries of a $2n \times 2n$ real symmetric positive definite matrix and its symplectic eigenvalues are known today. 
    These weak supermajorization relations are symplectic analogs of the classic Schur--Horn theorem
    described as follows.    
    Let $A \in \operatorname{Pd}(2n)$, and $d_{{\operatorname{s}}}(A)$ denote the $n$-vector whose entries are given by the symplectic eigenvalues $d_1(A),\ldots, d_n(A)$ of $A$.
    Express $A$ in the following block form:
    \begin{align}
    A = \begin{pmatrix}
        A_{11} & A_{12} \\
        A_{12}^T & A_{22}
    \end{pmatrix},\label{eq:block-form-A}
    \end{align}
    where the blocks $A_{11}, A_{12},$ and $A_{22}$ have size $n \times n$. 
    Let $\Delta_{11}, \Delta_{12},$ and $\Delta_{22}$ denote $n$-vectors consisting of the diagonal entries of $A_{11}, A_{12},$ and $A_{22}$, respectively.
    Define
    \begin{align}
    \Delta_{\operatorname{c}}(A) &\coloneqq \dfrac{\Delta_{11}+\Delta_{22}}{2} , \\
    \Delta_{\operatorname{s}}(A) &\coloneqq \sqrt{\Delta_{11} \cdot \Delta_{22}} , \\
    \Delta_{\operatorname{w}}(A) &\coloneqq \sqrt{\dfrac{\Delta_{11}^2  + \Delta_{22}^2}{2}}, \label{eq:definition_Delta_w} \\
    \Delta_{\operatorname{h}}(A) &\coloneqq \sqrt{\dfrac{\Delta_{11}^2  + \Delta_{22}^2 + 2 \Delta_{12}^2 }{2}}. \label{eq:definition_Delta_h}
    \end{align}
    Here the sums, products, divisions, squares, and square-roots of the vectors are taken entry-wise.
    Bhatia and Jain \cite{bhatia2020schur} proved the following weak supermajorization relations:
    \begin{align}
        \Delta_{\operatorname{c}}(A) &\prec^w d_{\operatorname{s}}(A), \label{eq:symplectic-schur-c}\\
        \Delta_{\operatorname{s}}(A) &\prec^w d_{\operatorname{s}}(A); \label{eq:symplectic-schur-s}
    \end{align}
     and Huang \cite{huang2023} showed another set of such relations:
     \begin{align}
        \Delta_{\operatorname{w}}(A) &\prec^w d_{\operatorname{s}}(A), \label{eq:symplectic-schur-w}\\
        \Delta_{\operatorname{h}}(A) &\prec^w d_{\operatorname{s}}(A).\label{eq:symplectic-schur-h}
    \end{align}
    
    It was shown in \cite{mishra2024equality} that majorization holds in the weak supermajorization \eqref{eq:symplectic-schur-c} if and only if $A$ is orthosymplectically diagonalizable in the sense of Williamson's theorem.
    In the following theorems, we establish precise conditions for majorization to hold in the remaining weak supermajorization relations \eqref{eq:symplectic-schur-s}, \eqref{eq:symplectic-schur-w}, and \eqref{eq:symplectic-schur-h}.
    \begin{theorem}\label{thm:majorization_delta_s_diagonal}
        Let $A \in \operatorname{Pd}(2n)$, and define
        \begin{align}
            M \coloneqq \operatorname{diag}\left(\sqrt[4]{\Delta_{11}^{-1} \cdot \Delta_{22}} \right) \oplus \operatorname{diag}\left(\sqrt[4]{\Delta_{11} \cdot \Delta_{22}^{-1}} \right) \in \operatorname{Sp}(2n).\label{eq:symplectic_matrix_constructed_from_diagonal_A}
        \end{align} 
        We have $\Delta_{\operatorname{s}}(A) \prec d_{\operatorname{s}}(A)$ if and only if $M A M$ is orthosymplectically diagonalizable in the sense of Williamson's theorem.
    \end{theorem}
    \begin{theorem}\label{thm:majorization_delta_w_diagonal}
        Let $A \in \operatorname{Pd}(2n)$.
        The following statements are equivalent:
        \begin{enumerate}
            \item [(i)] $\Delta_{\operatorname{w}}(A) \prec d_{\operatorname{s}}(A)$;
            \item [(ii)] $\Delta_{\operatorname{h}}(A) \prec d_{\operatorname{s}}(A)$;
            \item [(iii)] $A$ is orthosymplectically diagonalizable in the sense of Williamson's theorem.
        \end{enumerate}
    \end{theorem}
    \begin{remark}
        It is interesting to note that the necessary and sufficient conditions for both $\Delta_{\operatorname{w}}(A) \prec d_{\operatorname{s}}(A)$ and $\Delta_{\operatorname{h}}(A) \prec d_{\operatorname{s}}(A)$ are the same, which is also the same for the majorization condition $\Delta_{\operatorname{c}}(A) \prec d_{\operatorname{s}}(A)$ established in \cite{mishra2024equality}.
    \end{remark}
\subsection{Weak supermajorizations in symplectic pinching}
    Let $m_1,\ldots, m_k$ be positive integers and $n=m_1+\cdots+m_k$.
    Let $X_i$ be an $m_i \times m_i$ matrix for $i=1,\ldots, k$.
    Denote by $\oplus X_i$ the usual direct sum of the matrices $X_1,\ldots, X_k$.
    Suppose $A_i$ is a $2m_i \times 2m_i$ matrix partitioned into blocks as 
    \begin{align}
        A_i = 
            \begin{pmatrix}
                E_i & F_i \\
                G_i & H_i
            \end{pmatrix},
    \end{align}
    where $E_i, F_i, G_i, H_i$ are matrices of size $m_i \times m_i$ for $i=1,\ldots, k$.
    Define the $\operatorname{s}$-direct sum of $A_1,\ldots, A_k$ by
    \begin{align}
        \oplus^{\operatorname{s}} A_i 
            \coloneqq
                \begin{pmatrix}
                \oplus E_i & \oplus F_i \\
                \oplus G_i & \oplus H_i
            \end{pmatrix}.            
    \end{align}
    
    Suppose $X$ is an $n \times n$ matrix partitioned into blocks as $X=[X_{ij}]$, where each diagonal block $X_{ii}$ is an
    $m_i \times m_i$ matrix.
    A pinching of $X$ is defined as the block-diagonal matrix $\mathscr{C}\!(X) \coloneqq \oplus X_{ii}$.
    Let $A$ be a $2n \times 2n$ matrix in the block form
    \begin{align}
    A \coloneqq 
        \begin{pmatrix}
            E & F \\
            G & H
        \end{pmatrix},
    \end{align}
    where $E=\left[E_{ij} \right], F=\left[F_{ij} \right], G=\left[G_{ij} \right], H=\left[H_{ij} \right]$ are $n \times n$ block matrices with $m_i \times m_i$ diagonal blocks $E_{ii}, F_{ii}, G_{ii}, H_{ii}$, respectively.
    Define the $\operatorname{s}$-pinching of $A$ as
    \begin{align}\label{eq:symplectic_pinching_definition}
        \mathscr{C}^{\operatorname{s}}\!(A) \coloneqq 
        \begin{pmatrix}
            \mathscr{C}\!(E) & \mathscr{C}\!(F) \\
            \mathscr{C}\!(G) & \mathscr{C}\!(H)
        \end{pmatrix}.
    \end{align}
    We then have 
    \begin{align}\label{eq:symplectic_pinching_definition_s_direct_sum}
        \mathscr{C}^{\operatorname{s}}\!(A) = \oplus^{\operatorname{s}} A_i,
    \end{align}
    where $A_i \coloneqq \begin{psmallmatrix} E_{ii} & F_{ii} \\ G_{ii} & H_{ii} \end{psmallmatrix}$ for $i=1,\ldots, k$.
    
    It was proved by Bhatia and Jain \cite{bhatia2015symplectic} that the following weak supermajorization relation holds: for $A \in \operatorname{Pd}(2n)$, we have
    \begin{align}\label{eq:symplectic_pinching_weak_supermajorization}
        d_{\operatorname{s}}(\mathscr{C}^{\operatorname{s}}\!(A)) \prec^{w} d_{\operatorname{s}}(A).
    \end{align}
    In the following theorem, we provide a necessary and sufficient condition for the majorization to hold in the weak supermajorization \eqref{eq:symplectic_pinching_weak_supermajorization}.
    \begin{theorem}\label{thm:s_pinching_majorization_condition}
        Let $A \in \operatorname{Pd}(2n)$ be given. 
        Let $\mathscr{C}^{\operatorname{s}}\!(A)=\oplus^{\operatorname{s}} A_i$ be the symplectic pinching of $A$ as described in \eqref{eq:symplectic_pinching_definition_s_direct_sum}.
        Let $M_{i} \in \operatorname{Sp}(2m_i; A_i)$ be arbitrarily chosen for $i=1,\ldots, k$, and set
        $M \coloneqq \oplus^{\operatorname{s}} M_i$.
        We then have $d_{\operatorname{s}}(\mathscr{C}^{\operatorname{s}}\!(A)) \prec d_{\operatorname{s}}(A)$ if and only if $M^TAM$ is orthosymplectically diagonalizable in the sense of Williamson's theorem.
    \end{theorem}

    Let $A \in \operatorname{Pd}(2n)$ be partitioned into blocks given by \eqref{eq:block-form-A}. 
    Define the symplectic diagonal of $A$ as
    \begin{align}
        \mathfrak{D}^{\operatorname{s}}\!(A) 
            \coloneqq 
                \begin{pmatrix}
                    \operatorname{diag}(\Delta_{11}) & \operatorname{diag}(\Delta_{12}) \\
                    \operatorname{diag}(\Delta_{12}) & \operatorname{diag}(\Delta_{22})
                \end{pmatrix}.
    \end{align}
    Observe that $\mathfrak{D}^{\operatorname{s}}\!(A)$ is the $\operatorname{s}$-pinching of $A$ for $k=n$ whence $m_i=1$ for all $i=1,\ldots, n$.  
    As a direct consequence of Theorem~\ref{thm:s_pinching_majorization_condition}, we have
    \begin{corollary}\label{cor:sym_diag_majorization_condition}
        Let $A \in \operatorname{Pd}(2n)$ whose symplectic diagonal is given by $\mathfrak{D}^{\operatorname{s}}\!(A) = \oplus^{\operatorname{s}} A_i$, where $A_i$ are of the form 
        \begin{align}\label{eq:matrix_Ai}
            A_i = 
            \begin{pmatrix}
                \alpha_i & \beta_i \\
                \beta_i & \gamma_i
            \end{pmatrix},
        \end{align}
        $\alpha_i, \gamma_i > 0$ and $\alpha_i \gamma_i - \beta_i^2 >0$ for all $i=1,\ldots, n$.
        We have $M_i \in \operatorname{Sp}(2; A_i)$ given by
        \begin{align}\label{eq:sym_mat_diagonalizing_Ai}
            M_i 
                \coloneqq 
                \dfrac{1}{\sqrt[4]{\alpha_i \gamma_i(\alpha_i \gamma_i-\beta_i^2)}}
                \begin{pmatrix}
                    \sqrt[4]{ \gamma_i/\alpha_i} & 0 \\
                    0 & \sqrt[4]{\alpha_i/ \gamma_i}
                \end{pmatrix}
                \begin{pmatrix}
                    \sqrt{\alpha_i \gamma_i-\beta_i^2} & -\beta_i \\
                    0 & \sqrt{\alpha_i\gamma_i}
                \end{pmatrix}
        \end{align}
        for $i=1,\ldots, n$.
        Let $M \coloneqq \oplus^{\operatorname{s}} M_i$.
        Then $d_{\operatorname{s}}(\mathfrak{D}^{\operatorname{s}}\!(A)) \prec d_{\operatorname{s}}(A)$
        if and only if $M^T A M$ is orthosymplectically diagonalizable in the sense of Williamson's theorem.
    \end{corollary}
\section{Proofs}
    It will be helpful to observe the following.
    Let $A \in \operatorname{Pd}(2n)$ be fixed.
    By Williamson's theorem, there exists a symplectic matrix $M \in \operatorname{Sp}(2n)$ such that $A=M^{-T}(D \oplus D)(M^{-T})^T$, where $D$ is an $n \times n$ diagonal matrix consisting of the symplectic eigenvalues of $A$. 
    Set $N \coloneqq M^{-T} \in \operatorname{Sp}(2n)$. 
    Write this matrix in the block form
    \begin{align}\label{eq:N_block_form_symplectic}
    N =   
        \begin{pmatrix}
            	P & Q \\
            	R & S
        \end{pmatrix},
    \end{align}
    where $P,Q,R,S$ are $n \times n$ matrices.  
    We thus get
    \begin{align}
    A = N(D \oplus D)N^T  
        &= 
        \begin{pmatrix}
            PDP^T+ QDQ^T & PDR^T+QDS^T \\
            RDP^T+SDQ^T & RDR^T + SDS^T
        \end{pmatrix}.
    \end{align}
    Recall from the block-form \eqref{eq:block-form-A} of $A$ that
    \begin{align}
        A_{11} &= PDP^T+ QDQ^T,\\
        A_{12} &=PDR^T+QDS^T, \\
        A_{22} &=RDR^T + SDS^T.
    \end{align}
    The diagonal vectors of these matrices are given by
    \begin{align}
        \Delta_{11} &= (P\circ P + Q \circ Q ) d_{\operatorname{s}}(A), \label{eq:Delta_one_one_representation} \\
        \Delta_{12} &= (P\circ R + Q \circ S ) d_{\operatorname{s}}(A), \label{eq:Delta_one_two_representation} \\
        \Delta_{22} &= (R\circ R + S \circ S ) d_{\operatorname{s}}(A), \label{eq:Delta_two_two_representation}
    \end{align}
    where $\circ$ denotes the Hadamard (entry-wise) product of matrices. 
    The symplectic matrix $N$ is orthogonal if and only if $P=S$, $Q=-R$, and $P+\iota Q$ is an $n \times n$ complex unitary matrix, where $\iota \coloneqq \sqrt{-1}$ is the imaginary unit.
    This is stated in \cite[Section~5]{bhatia2015symplectic}. 
    See also \cite[Proposition~2.12]{degosson} for a proof.
\subsection{Proof of Theorem~\ref{thm:majorization_delta_s_diagonal}}
    One can verify that the matrix $M$ defined in \eqref{eq:symplectic_matrix_constructed_from_diagonal_A} is a symplectic diagonal matrix.
    Thus, we have $d_{\operatorname{s}}(A)=d_{\operatorname{s}}(M^TAM)=d_{\operatorname{s}}(MAM)$ because symplectic eigenvalues are invariant under the symplectic orbit \cite{bhatia2015symplectic}.
    Also, by construction, we have
    \begin{align}
        \Delta_{\operatorname{s}}(A) &= \Delta_{\operatorname{c}}(MAM).
    \end{align}
    Therefore, the majorization condition $\Delta_{\operatorname{s}}(A) \prec d_{\operatorname{s}}(A)$ is equivalent to $\Delta_{\operatorname{c}}(MAM) \prec d_{\operatorname{s}}(MAM)$.
    We know from \cite{mishra2024equality} that $\Delta_{\operatorname{c}}(MAM) \prec d_{\operatorname{s}}(MAM)$ if and only if $MAM$ is orthosymplectically diagonalizable in the sense of Williamson's theorem.
    This completes the proof.

\subsection{Proof of Theorem~\ref{thm:majorization_delta_w_diagonal}}
    We first prove that $(ii) \Rightarrow (i)$.
    Let us assume that the majorization $\Delta_{\operatorname{h}}(A) \prec d_{\operatorname{s}}(A)$ holds.
    By definition, we have $\Delta_{\operatorname{w}}(A) \leq \Delta_{\operatorname{h}}(A)$.
    It then follows by \eqref{eq:symplectic-schur-w} that $\Delta_{\operatorname{w}}(A) \prec d_{\operatorname{s}}(A)$.

    We now show that $(i)\Rightarrow (iii)$.
    Assume that the majorization $\Delta_{\operatorname{w}}(A) \prec d_{\operatorname{s}}(A)$ holds.
    One can easily verify that $\Delta_{\operatorname{c}}(A) \leq \Delta_{\operatorname{w}}(A)$.
    It then follows by \eqref{eq:symplectic-schur-c} that $\Delta_{\operatorname{c}}(A) \prec d_{\operatorname{s}}(A)$. This implies from \cite{mishra2024equality} that $A$ is orthosymplectically diagonalizable.

    At last, we establish $(iii) \Rightarrow (ii)$.
    Assume that $A$ is diagonalizable by an orthosymplectic matrix in the sense of Williamson's theorem.
    Let $N$ be an orthosymplectic matrix given in the block form \eqref{eq:N_block_form_symplectic}, and hence $P=S$ and $Q=-R$.
    From \eqref{eq:Delta_one_one_representation}, \eqref{eq:Delta_one_two_representation}, and \eqref{eq:Delta_two_two_representation} we thus get $\Delta_{12}=0$ and
    \begin{align}
        \Delta_{11}=\Delta_{22} &= \left(P\circ P  + Q\circ Q \right)d_{\operatorname{s}}(A).
    \end{align}
    Then, \eqref{eq:definition_Delta_h} implies
    \begin{align}
        \Delta_{\operatorname{h}}(A) = \left(P\circ P  + Q\circ Q \right)d_{\operatorname{s}}(A).\label{eq:delta_h_dstochastic_symplectic_vector}
    \end{align}    
    Since $N$ is an orthosymplectic matrix, the matrix $P\circ P  + Q\circ Q$ is doubly stochastic \cite{bhatia2015symplectic}.
    So, \eqref{eq:delta_h_dstochastic_symplectic_vector} implies that $\Delta_{\operatorname{h}}(A) \prec d_{\operatorname{s}}(A)$.
    This concludes the proof.
    
\subsection{Proof of Theorem~\ref{thm:s_pinching_majorization_condition}}
    Let $\mathscr{C}^{\operatorname{s}}\!(A)=\oplus^{\operatorname{s}} A_i$ be the $\operatorname{s}$-pinching of $A$, and $M= \oplus^{\operatorname{s}} M_i$ as given in the theorem.
    The symplectic eigenvalues of $\mathscr{C}^{\operatorname{s}}\!(A)$ are the symplectic eigenvalues of $A_1,\ldots, A_k$ put together. 
    Moreover, we have that $M \in \operatorname{Sp}(2n; \mathscr{C}^{\operatorname{s}}\!(A))$.
    See \cite[Section~6]{bhatia2015symplectic}.
    
    Assume that $d_{\operatorname{s}}(\mathscr{C}^{\operatorname{s}}\!(A)) \prec d_{\operatorname{s}}(A)$.
    We then have
    \begin{align}
        \operatorname{Tr}(M^T A M) 
            &= \sum_{i=1}^k \operatorname{Tr}(M_i^T A_i M_i) \\
            &= \operatorname{Tr}(M^T \mathscr{C}^{\operatorname{s}}\!(A) M) \\
            &= 2\sum_{i=1}^n d_i(\mathscr{C}^{\operatorname{s}}\!(A)).\label{eq:trace_m_transpose_A_m_symplectic_pinching}
    \end{align}
    The assumption $d_{\operatorname{s}}(\mathscr{C}^{\operatorname{s}}\!(A)) \prec d_{\operatorname{s}}(A)$ thus implies that
    \begin{align}
        \operatorname{Tr}(M^T A M) = 2\sum_{i=1}^n d_i(A).\label{eq:trace_m_transpose_A_m_symplectic}
    \end{align}
    By Theorem~4.6~(ii) of \cite{sags_2021}, the matrix $M^T A M$ is orthosymplectically diagonalizable in the sense of Williamson's theorem.

    Converse is rather straightforward.
    Indeed, if the matrix $M^T A M$ is orthosymplectically diagonalizable in the sense of Williamson's theorem, then \eqref{eq:trace_m_transpose_A_m_symplectic} holds.
    Combining \eqref{eq:trace_m_transpose_A_m_symplectic} with \eqref{eq:trace_m_transpose_A_m_symplectic_pinching} then gives that the sums of the symplectic eigenvalues of $A$ and $\mathscr{C}^{\operatorname{s}}\!(A)$ are equal.
    This combined with \eqref{eq:symplectic_pinching_weak_supermajorization} proves that $d_{\operatorname{s}}(\mathscr{C}^{\operatorname{s}}\!(A)) \prec d_{\operatorname{s}}(A)$.

\subsection{Proof of Corollary~\ref{cor:sym_diag_majorization_condition}}
        Set 
        \begin{align}
            L_i 
                &\coloneqq 
                \begin{pmatrix}
                    \sqrt[4]{ \gamma_i/\alpha_i} & 0 \\
                    0 & \sqrt[4]{\alpha_i/ \gamma_i}
                \end{pmatrix}, \\
            N_i
                &\coloneqq  
                    \dfrac{1}{\sqrt[4]{\alpha_i \gamma_i(\alpha_i \gamma_i-\beta_i^2)}}
                    \begin{pmatrix}
                    \sqrt{\alpha_i \gamma_i-\beta_i^2} & -\beta_i \\
                    0 & \sqrt{\alpha_i\gamma_i} 
                    \end{pmatrix},
        \end{align}
        so that $M_i = L_i N_i$.
        It is easy to verify that a $2 \times 2$ real matrix is symplectic if and only if its determinant is equal to one.
        This implies that $L_i, N_i$ are symplectic matrices, and hence $M_i$ is a symplectic matrix. 
        The fact that $M_i$ diagonalizes the matrix $A_i$ given in \eqref{eq:matrix_Ai} in the sense of Williamson's theorem can be seen as follows.
        We have
        \begin{align}
            M_i^T A_i M_i 
                &= N_i^T L_i^T A_i L_i N_i \\
                &= N_i^T \begin{pmatrix}
                    \sqrt[4]{ \gamma_i/\alpha_i} & 0 \\
                    0 & \sqrt[4]{\alpha_i/ \gamma_i}
                    \end{pmatrix} \begin{pmatrix}
                    \alpha_i & \beta_i \\
                    \beta_i & \gamma_i
                    \end{pmatrix} \begin{pmatrix}
                    \sqrt[4]{ \gamma_i/\alpha_i} & 0 \\
                    0 & \sqrt[4]{\alpha_i/ \gamma_i}
                    \end{pmatrix} N_i \\
                &= N_i^T \begin{pmatrix}
                    \sqrt{\alpha_i \gamma_i} & \beta_i \\
                    \beta_i & \sqrt{\alpha_i \gamma_i}
                    \end{pmatrix} N_i \\
                &= \dfrac{1}{\sqrt{\alpha_i \gamma_i(\alpha_i \gamma_i-\beta_i^2)}}
                \begin{pmatrix}
                    \sqrt{\alpha_i \gamma_i-\beta_i^2} & 0 \\
                    -\beta_i & \sqrt{\alpha_i\gamma_i} 
                \end{pmatrix} \nonumber \\
                & \hspace{3cm}
                \begin{pmatrix}
                \sqrt{\alpha_i \gamma_i} & \beta_i \\
                \beta_i & \sqrt{\alpha_i \gamma_i}
                \end{pmatrix}
                \begin{pmatrix}
                    \sqrt{\alpha_i \gamma_i-\beta_i^2} & -\beta_i \\
                    0 & \sqrt{\alpha_i\gamma_i} 
                \end{pmatrix} \\
                &= \dfrac{1}{\sqrt{\alpha_i \gamma_i(\alpha_i \gamma_i-\beta_i^2)}}
                    \begin{pmatrix}
                    \sqrt{\alpha_i \gamma_i} \left(\alpha_i \gamma_i-\beta_i^2 \right) & 0 \\
                    0 & \sqrt{\alpha_i \gamma_i} \left(\alpha_i \gamma_i-\beta_i^2 \right)
                \end{pmatrix} \\
                &= \begin{pmatrix}
                    \sqrt{\alpha_i \gamma_i-\beta_i^2 } & 0 \\
                    0 & \sqrt{\alpha_i \gamma_i-\beta_i^2 }
                \end{pmatrix}.
        \end{align}
    We have thus shown that $M_i \in \operatorname{Sp}(2; A_i)$ for all $i=1,\ldots, n$.
    The assertion of the corollary now follows directly from Theorem~\ref{thm:s_pinching_majorization_condition}.





\section{Conclusion}
    We established necessary and sufficient conditions for the majorization to hold in four types of weak supermajorization relations of symplectic eigenvalues: three weak supermajorizations given by symplectic Schur--Horn theorems \eqref{eq:symplectic-schur-s}, \eqref{eq:symplectic-schur-w}, \eqref{eq:symplectic-schur-h} and one weak supermajorization corresponding to $\operatorname{s}$-pinching \eqref{eq:symplectic_pinching_weak_supermajorization}.





\section*{Acknowledgments}
    The authors thank an anonymous reviewer for pointing out some typos in the first draft of the paper, and for their helpful suggestions that increased readability of the paper.
    The work of Shaowu Huang is supported by the National Natural Science Foundation of China (Grant No. 12201332). Hemant K. Mishra acknowledges supports from the NSF under grant no.~2304816 and AFRL under agreement no.~FA8750-23-2-0031.





\begin{Backmatter}

\bibliographystyle{unsrt}
\bibliography{RefBT}

\printaddress
\end{Backmatter}

\end{document}